# Gradient algorithms for finding common Lyapunov functions

Daniel Liberzon*

Coordinated Science Laboratory
Univ. of Illinois at Urbana-Champaign
Urbana, IL, U.S.A.

Roberto Tempo†

IRITI-CNR
Politecnico di Torino
Torino, Italy

October 25, 2018

**Abstract**

This paper is concerned with the problem of finding a quadratic common Lyapunov function for a family of stable linear systems. We present gradient iteration algorithms which give deterministic convergence for finite system families and probabilistic convergence for infinite families.

## 1 Introduction

A switched system is a dynamical system described by a family of continuous-time subsystems and a rule that governs the switching between them. Such systems arise as models of processes regulated by switching control mechanisms and/or affected by abrupt changes in the dynamics. It is well known and easy to demonstrate that switching between stable subsystems may lead to instability. This makes stability analysis of switched systems an important and challenging problem, which has received considerable attention in the recent literature. We refer the reader to the survey articles [12] and [5] for an overview of available results.

Stability of a switched system is guaranteed if the individual subsystems share a common Lyapunov function. If the system enjoys this property, then stability is preserved for arbitrary switching sequences. When the subsystems being switched are obtained as feedback interconnections of a given process with different stabilizing controllers, this means that one does not need to worry about stability and can concentrate on other issues such as performance.

A particular case of interest is when the subsystems are linear time-invariant and a quadratic common Lyapunov function is sought. Although a number of conditions for the existence of such a Lyapunov function have been obtained, general results are lacking. On the other hand, the problem of finding a quadratic common Lyapunov function amounts to solving a system of linear matrix inequalities (LMIs), and efficient methods for solving such inequalities are available [2]. However, this approach becomes infeasible when the number of subsystems being switched is large, and is not useful in the case of an infinite family of subsystems.

---
*Email: liberzon@uiuc.edu. Supported by NSF Grants ECS-0134115 and ECS-0114725.
†Email: tempo@polito.it.

In view of the above remarks, there is a need for developing computationally tractable algorithms which can be used to find quadratic common Lyapunov functions for large families of linear systems. It turns out that this can be achieved by handling matrix inequality constraints sequentially rather than simultaneously. (For comparison, we mention the iterative algorithm proposed in [15] for a special case of the above problem, namely, finding a common Lyapunov function for a finite family of commuting linear systems.) In the case of an infinite family of systems, one can relax the objective and apply the same ideas combined with randomization to obtain probabilistic convergence results.

The use of gradient algorithms for solving matrix inequalities was proposed in [3] and [18] in the context of probabilistic control design for uncertain linear systems; see also [6, 8, 16, 19] for related recent developments. The main ideas go back to the early work on solving algebraic inequalities reported in [1, 14, 17, 20]. The goal of this paper is to show how the problem posed above can be addressed using similar techniques. We describe a general framework for developing gradient iteration algorithms for finding quadratic common Lyapunov functions and discuss and compare several specific design choices, some of which were used in previous work cited above while others were not. Convergence to a quadratic common Lyapunov function (if one exists) in a finite number of steps is guaranteed for a finite family of linear systems and takes place with probability one for an infinite family. While our approach parallels that of [18], there are several important technical differences on which we elaborate below. We first treat finite families and then discuss an extension to infinite families.

## 2   Problem formulation and notation

Suppose first that we are given a family of real Hurwitz $n \times n$ matrices $A_1, \ldots, A_N$, where $n$ and $N$ are positive integers. We write $P \geq 0$ (or $P > 0$) to indicate that a matrix $P$ is symmetric nonnegative definite (respectively, positive definite), and $P \leq 0$ (or $P < 0$) to indicate that $P$ is symmetric nonpositive definite (respectively, negative definite). We assume throughout that there exists a matrix $P > 0$ which satisfies

$$PA_i + A_i^T P < 0, \qquad i = 1, \ldots, N.$$

This means that the quadratic function $V(x) := x^T P x$ is a common Lyapunov function for the family of asymptotically stable linear systems

$$\dot{x} = A_i x, \qquad i = 1, \ldots, N. \tag{1}$$

Fix an arbitrary matrix $Q > 0$. Multiplying $P$ by a sufficiently large positive number, we see that the system of inequalities

$$PA_i + A_i^T P + Q \leq 0, \qquad i = 1, \ldots, N \tag{2}$$

has a solution $P > 0$. Moreover, if a symmetric matrix $P$ satisfies (at least one of) the inequalities (2), then it is well known that we automatically have $P > 0$; see, e.g., [4, p. 132]. Thus the problem of finding a quadratic common Lyapunov function for the family (1) is equivalent to that of finding a symmetric matrix $P$ which satisfies (2). In what follows, we will be concerned with the latter problem. We denote the set of symmetric solutions of the inequalities (2) by $\mathcal{L}$.

The space of real symmetric $n \times n$ matrices is equipped with the Frobenius norm $||R|| := (\sum_{i,j=1}^n R_{ij}^2)^{1/2}$ and the inner product $\langle R, S \rangle := \text{tr } RS$. It is important to note that since the set $\mathcal{L}$ is nonempty by assumption, it must have a nonempty interior. Indeed, if $P \in \mathcal{L}$, then a standard perturbation argument

such as the one given in [9, Example 5.1] can be used to show that $\mathcal{L}$ contains a neighborhood of $\gamma P$ for every $\gamma > 1$. In fact, we have $\gamma P + \Delta P \in \mathcal{L}$ for every symmetric matrix $\Delta P$ satisfying the bound

$$\lambda_{\max}(\Delta P) \leq \frac{(\gamma - 1)\lambda_{\min}(Q)}{2\max_{i=1,\ldots,N} \sigma_{\max}(A_i)}$$

where $\lambda_{\max}$ and $\lambda_{\min}$ denote the largest and the smallest eigenvalue, respectively, and $\sigma_{\max}$ denotes the largest singular value. It is well known that $\lambda_{\max}(\Delta P) \leq \|\Delta P\|$; see, e.g., [7, pp. 296–297]. Since $\gamma$ can be arbitrarily large, we see that for every $r > 0$ there exists a ball of radius $r$ which is contained in $\mathcal{L}$. This is in contrast with the case of Riccati inequalities treated in [18], where a sufficiently small $r$ with the above property is assumed to exist but is not explicitly known.

We will need the notion of projection of a symmetric matrix $R$ onto the convex cone of nonnegative definite matrices. This projection is defined as

$$R^+ := \arg\min_{S \geq 0} \|R - S\|.$$

It can be computed explicitly as follows (see [17]): If $R = U\Lambda U^T$, where $U$ is orthogonal and $\Lambda$ is diagonal with entries $\lambda_1, \ldots, \lambda_n$, then $R^+ = U\Lambda^+ U^T$, where $\Lambda^+$ is diagonal with entries $\max\{0, \lambda_1\}, \ldots, \max\{0, \lambda_n\}$. We also denote $R - R^+$ by $R^-$.

## 3 Gradient algorithms

Let $f$ be a convex functional on the space of symmetric matrices, which assigns to a matrix $R$ a real number $f(R)$, with the property that $f(R) \leq 0$ if and only if $R \leq 0$. We suppose that this functional is differentiable (this condition can be relaxed somewhat, as discussed later) and that its gradient is symmetric. We denote this gradient by $\partial_R f$, and similar notation will be used for other gradients appearing below. Specific examples of functionals with the above properties will be given in Section 5; here we keep the discussion general.

Given a symmetric matrix $P$ and another matrix $A$, we let

$$v(P, A) := f(PA + A^T P + Q)$$

where $f$ is the functional introduced above and $Q > 0$ is the matrix from (2). Since $f$ is convex, $v$ is convex in $P$. Well-known results imply that for each integer $i$ between 1 and $N$, solutions of the gradient system $\dot{P} = -\partial_P v(P, A_i)$ converge to the set $\{P : v(P, A_i) \leq 0\} = \{P : PA_i + A_i^T P + Q \leq 0\}$. The same is true for the associated discrete iterations $P_{k+1} = P_k - \mu_k \partial_P v(P_k, A_i)$, $k = 0, 1, \ldots$ if the step-sizes $\mu_k$ are chosen appropriately. Moreover, by switching between the above iterations for different values of $i$, we can make $P_k$ converge to the intersection of the corresponding sets, which is precisely the set $\mathcal{L}$ of solutions of (2). This happens because the distance from $P_k$ to $\mathcal{L}$ with respect to the Frobenius norm is a decreasing function of $k$.

We now make the above discussion precise by describing how the gradient iterations are to be carried out to ensure convergence to $\mathcal{L}$ in a finite number of steps. We need to pick a "scheduling function" $h$ from nonnegative integers to the set $\{1, \ldots, N\}$ which has the following *revisitation property*: For every integer $i$ between 1 and $N$ and for every integer $l \geq 0$ there exists an integer $k \geq l$ such that $h(k) = i$. An obvious choice is $h(k) := (k \mod N) + 1 = k - N\lfloor k/N \rfloor + 1$.

For each $k = 0, 1, \ldots$, let us define the step-size by the formula

$$\mu_k := \frac{\alpha v(P_k, A_{h(k)}) + r\|\partial_P v(P_k, A_{h(k)})\|}{\|\partial_P v(P_k, A_{h(k)})\|^2} \tag{3}$$

where $0 \leq \alpha \leq 1$ and $r > 0$ are arbitrary. Consider the iterations

$$P_{k+1} = \begin{cases} P_k - \mu_k \partial_P v(P_k, A_{h(k)}) & \text{if } v(P_k, A_{h(k)}) > 0, \\ P_k & \text{otherwise.} \end{cases} \tag{4}$$

We take the initial condition $P_0$ to be symmetric. (For example, a solution of one of the inequalities (2) provides a convenient choice for $P_0$.) Then $P_k$ is symmetric for each $k$, since $\partial_P v$ is symmetric in light of the following lemma and the fact that $\partial_R f$ is symmetric.

**Lemma 1** *The gradient of $v$ is given by*

$$\partial_P v(P, A) = A \partial_R f(PA + A^T P + Q) + \partial_R f(PA + A^T P + Q) A^T. \tag{5}$$

PROOF. Denoting by $\Delta P$ a small perturbation in $P$ and by $\approx$ equality up to first-order terms in $\Delta P$, we write

$$\begin{aligned} v(P + \Delta P, A) &= f(PA + A^T P + Q + \Delta PA + A^T \Delta P) \\ &\approx f(PA + A^T P + Q) + \langle \partial_R f(PA + A^T P + Q), \Delta PA + A^T \Delta P \rangle \\ &= v(P, A) + \operatorname{tr}\left( \left(A \partial_R f(PA + A^T P + Q) + \partial_R f(PA + A^T P + Q) A^T\right) \Delta P \right) \end{aligned}$$

from which (5) follows. $\square$

On the other hand, $P_{k+1}$ is not guaranteed to be positive definite or at least nonnegative definite, even if $P_k > 0$. To make sure that a nonnegative definite matrix is generated at every step, we could instead consider

$$P_{k+1} = \begin{cases} [P_k - \mu_k \partial_P v(P_k, A_{h(k)})]^+ & \text{if } v(P_k, A_{h(k)}) > 0, \\ P_k & \text{otherwise} \end{cases} \tag{6}$$

using the projection operation defined in Section 2. Since all matrices in the set $\mathcal{L}$ are positive definite, this modification also has the potential of improving convergence; see the proof of Theorem 1 below.

**Remark 1** The algorithm (6) exactly parallels the one proposed in [18] for finding nonnegative definite solutions of systems of Riccati inequalities. The Lyapunov inequalities (2), on the other hand, have the property that if a symmetric matrix $P_k$ satisfies at least one of them, then necessarily $P_k \geq 0$. This means that the projection is not really needed, and the convergence result presented in the next section implies that the algorithm (4) generates only nonnegative definite matrices after sufficiently many steps.

## 4 Deterministic convergence for finite families

We now demonstrate that the above gradient algorithms indeed provide convergence to the desired set $\mathcal{L}$ in a finite number of steps. When for a given $k$ we have $v(P_k, A_{h(k)}) > 0$, we say that a *correction step* is executed.

**Theorem 1** *When the algorithm* (4) *or the algorithm* (6) *is applied with the step-size given by* (3), *there exists an integer $k^*$ such that $P_{k^*} \in \mathcal{L}$.*

PROOF. The proof is carried out along the lines of [18]. Consider a $k$ for which $v(P_k, A_{h(k)}) > 0$ and so a correction step is executed. As shown in Section 2, the set $\mathcal{L}$ contains a ball of radius $r$, centered at some matrix $P^*$. We will prove that

$$\|P_{k+1} - P^*\|^2 \leq \|P_k - P^*\|^2 - r^2. \tag{7}$$

Since the revisitation property of $h$ guarantees that a correction step occurs at least once in every $N$ steps until $P_k \in \mathcal{L}$, we can conclude that no more than $N \lfloor \|P_0 - P^*\|^2 / r^2 \rfloor$ steps are needed, and the proof will be complete.

Consider the algorithm (6). We have

$$\|P_{k+1} - P^*\|^2 = \|[P_k - \mu_k \partial_P v(P_k, A_{h(k)})]^+ - P^*\|^2 \leq \|P_k - \mu_k \partial_P v(P_k, A_{h(k)}) - P^*\|^2$$

where the last inequality follows from the definition of the projection and the supporting hyperplane theorem. Thus we see that it is enough to show (7) for the algorithm (4).

To this end, define

$$\bar{P} := P^* + \frac{r}{\|\partial_P v(P_k, A_{h(k)})\|} \partial_P v(P_k, A_{h(k)}) \in \mathcal{L}.$$

Then we write

$$\|P_{k+1} - P^*\|^2 = \|P_k - P^* - \mu_k \partial_P v(P_k, A_{h(k)})\|^2 = \|P_k - P^*\|^2 + \mu_k^2 \|\partial_P v(P_k, A_{h(k)})\|^2$$
$$- 2\mu_k \langle \partial_P v(P_k, A_{h(k)}), P_k - \bar{P} \rangle - 2\mu_k \langle \partial_P v(P_k, A_{h(k)}), \bar{P} - P^* \rangle.$$

We now consider the last two terms. Due to convexity of $v$ in $P$, we have

$$\langle \partial_P v(P_k, A_{h(k)}), P_k - \bar{P} \rangle \geq v(P_k, A_{h(k)}) \geq \alpha v(P_k, A_{h(k)})$$

while the definition of $\bar{P}$ gives

$$\langle \partial_P v(P_k, A_{h(k)}), \bar{P} - P^* \rangle = r \|\partial_P v(P_k, A_{h(k)})\|.$$

Therefore,

$$\|P_{k+1} - P^*\|^2 \leq \|P_k - P^*\|^2 + \mu_k^2 \|\partial_P v(P_k, A_{h(k)})\|^2 - 2\mu_k \Big(\alpha v(P_k, A_{h(k)}) + r \|\partial_P v(P_k, A_{h(k)})\|\Big).$$

Substituting the value of $\mu_k$ defined in (3), we obtain

$$\|P_{k+1} - P^*\|^2 \leq \|P_k - P^*\|^2 - \frac{\big(\alpha v(P_k, A_{h(k)}) + r \|\partial_P v(P_k, A_{h(k)})\|\big)^2}{\|\partial_P v(P_k, A_{h(k)})\|^2} \leq \|P_k - P^*\|^2 - r^2$$

and so the inequality (7) holds as claimed. □

Since the matrix $P^*$ used in the above proof is not known, the number $k^*$ may be difficult to estimate in practice. One can of course let the algorithm run for some number of steps $k$ and then check whether or not the matrix $P_k$ satisfies the inequalities (2); performing such a check is an easy task compared with that of solving the inequalities directly. We also have a lot of freedom in the choice of the step-size; more information on how to choose the step-size efficiently in gradient algorithms of this kind can be found, e.g., in [10, 17, 20].

# 5 Possible choices for $f$

We suggest two admissible choices of the functional $f$. One is given by

$$f(R) := \|R^+\|^2. \tag{8}$$

This functional has all the properties required in Section 3, as shown in [17]. For completeness, we sketch the argument below.

**Lemma 2** *The functional (8) is convex and differentiable, with gradient given by*

$$\partial_R f(R) = 2R^+. \tag{9}$$

PROOF. The calculations that follow rely on some standard properties of the projection and the cone of symmetric nonnegative definite matrices; see [17] for details. We write

$$f(R + \Delta R) = \|R + \Delta R - (R + \Delta R)^-\|^2 = \inf_{S \leq 0} \|R + \Delta R - S\|^2$$
$$\leq \|R + \Delta R - R^-\|^2 = \|R^+ + \Delta R\|^2 = \|R^+\|^2 + \langle 2R^+, \Delta R \rangle + \|\Delta R\|^2$$
$$\approx f(R) + \langle 2R^+, \Delta R \rangle.$$

On the other hand,

$$f(R + \Delta R) = \|R + \Delta R - (R + \Delta R)^-\|^2 = \|R^+ + R^- + \Delta R - (R + \Delta R)^-\|^2$$
$$= \|R^+\|^2 + \langle 2R^+, \Delta R \rangle + 2\langle R^+, R^- \rangle + 2\langle R^+, -(R + \Delta R)^- \rangle + \|R^- + \Delta R - (R + \Delta R)^-\|^2$$
$$\geq f(R) + \langle 2R^+, \Delta R \rangle.$$

This proves (9), and the inequality $f(R + \Delta R) \geq f(R) + \langle \partial_R f(R), \Delta R \rangle$ implies that $f$ is convex. □

Combining the formulas (5) and (9), we see that with this choice of $f$ the gradient of $v$ is given by

$$\partial_P v(P, A) = 2A(PA + A^T P + Q)^+ + 2(PA + A^T P + Q)^+ A^T.$$

This is the quantity that needs to be calculated at every step when implementing the algorithms of Section 3. As explained in Section 2, this is a routine calculation based on eigenvalue decomposition. An analogous construction was used in [18].

Another option is to let $f(R)$ be the largest eigenvalue of $R$:

$$f(R) := \lambda_{\max}(R). \tag{10}$$

The following result is standard; see, e.g, [7, p. 372].

**Lemma 3** *The functional (10) is convex and, when $\lambda_{\max}(R)$ is a simple eigenvalue, differentiable with gradient given by*

$$\partial_R f(R) = xx^T$$

*where $x$ is a unit eigenvector of $R$ with eigenvalue $\lambda_{\max}(R)$.*

With this choice of $f$, the gradient of $v$ is given by

$$\partial_P v(P, A) = Axx^T + xx^T A^T$$

and the projection is no longer required for implementation of the algorithm (4). A disadvantage of this $f$ is that it is not differentiable when $\lambda_{\max}(R)$ is not a simple eigenvalue. However, since a generic matrix has distinct eigenvalues, this problem can always be fixed by slightly perturbing the matrix $P_k$ if necessary.

# 6 Extension: probabilistic convergence for infinite families

Let us now consider the more general situation where we are given a compact set of real Hurwitz matrices $\mathcal{A} := \{A_p : p \in \mathcal{P}\}$, parameterized by some index set $\mathcal{P}$ which is in general infinite. It is still true that the problem of finding a matrix $P > 0$ which satisfies the inequalities

$$PA_p + A_p^T P < 0 \qquad \forall p \in \mathcal{P}$$

is equivalent to the problem of finding a symmetric matrix $P$ which satisfies the inequalities

$$PA_p + A_p^T P + Q \leq 0 \qquad \forall p \in \mathcal{P}$$

where $Q > 0$ is arbitrary. As before, we denote the set of such matrices $P$ by $\mathcal{L}$ and assume that it is nonempty; this again implies that $\mathcal{L}$ contains balls of arbitrarily large radii.

When the set of matrices is infinite, choosing a matrix at each step using a function $h$ with the revisitation property as in Section 3 is no longer possible. Instead, we use randomization. Namely, for each $k$ we randomly pick a matrix in $\mathcal{A}$ according to some probability distribution on $\mathcal{A}$ with the following property: Every subset of $\mathcal{A}$ which is open relative to $\mathcal{A}$ has a nonzero probability measure. We let $h(k)$ be the index of the matrix chosen in this way (this index may not be unique if the map $p \mapsto A_p$ is not injective). Note that $\mathcal{A}$ may be composed of several disjoint subsets such as intervals or isolated points, and we must assign a positive probability measure to each of them. If $\mathcal{A}$ consists of a finite number of isolated points, then randomization provides an alternative to the approach described in Section 3. As is well known, for infinite families described by convex polyhedra the problem reduces to the corresponding problem for finite families generated by the vertices; see, e.g., [13] and the references therein.

The randomized versions of the algorithms from Section 3 can now be defined by the same formulas (4) and (6) as before. Remark 1 also applies here. If $P_k \notin \mathcal{L}$, then $v(P_k, \bar{A}) > 0$ for some $\bar{A} \in \mathcal{A}$, and by continuity the same is true for all $A$ in a sufficiently small neighborhood of $\bar{A}$ in $\mathcal{A}$. Since the probability measure of this neighborhood is positive, a correction step will be executed with probability one after a finite number of steps. Proceeding exactly as in the proof of Theorem 1, we arrive at the following probabilistic counterpart of that theorem.

**Theorem 2** *When the algorithm* (4) *or the algorithm* (6) *is applied with the step-size given by* (3), *the probability of obtaining $P_{k^*} \in \mathcal{L}$ for some integer $k^*$ is equal to one.*

The above result parallels Theorem 1 from [18], although the context here is different. In [18] it is also shown how one can compute a lower bound on the probability of finding a solution in a given number of steps, provided that additional a priori information is available.

# 7 Concluding remarks

We presented gradient iteration algorithms for finding a quadratic common Lyapunov function for a family of asymptotically stable linear systems. We derived a deterministic convergence result for finite families, and then used randomization to obtain a probabilistic counterpart for infinite families. The algorithms were described on a general level with some freedom left in the design choices, several possibilities for which were also discussed.

We implemented the above algorithms in MATLAB for interval families of $4 \times 4$ and $5 \times 5$ upper-triangular Hurwitz matrices, for which quadratic common Lyapunov functions are known to exist (see, e.g., [11]). We applied the method of Section 3, noting that it suffices to work with the vertices. This was a straightforward task, and in all cases we observed convergence in a reasonable number of iterations: approximately 5,000 for the $4\times 4$ case and 75,000 for the $5\times 5$ case. The number of actual correction steps was an order of magnitude lower than the total number of iterations. Trade-offs among different choices of the specific gradient algorithm and the step-size remain to be understood. It is important to note that no analytical results on the existence of a common Lyapunov function are available for generic matrices of dimension higher than $2\times 2$, and that interval families of $4\times 4$ and $5\times 5$ triangular matrices give rise to systems of $2^{10}$ and $2^{15}$ LMIs, respectively, solving which simultaneously is computationally intractable.

The results reported here are valid under the assumption that a quadratic common Lyapunov function exists. They do not offer insight into the question of how to determine whether or not such a Lyapunov function can be found; this issue needs further investigation. By combining the above algorithms with a suitable recursive partitioning procedure, we can compute common Lyapunov functions for subsets of the given family of systems when the overall common Lyapunov function may not exist or is not found in a prescribed number of steps. Stability of the switched system is then guaranteed under appropriate constraints on the switching rate between different subsets; cf. [12]. Problems for future work also include choosing an optimal schedule for the iterations with respect to some performance criterion (for example, the speed of convergence) and incorporating additional constraints on the Lyapunov function to be found (such as the average or worst-case decay rate).

**Acknowledgment.** We are thankful to Chris Hadjicostis and Boris Polyak for helpful discussions. This research was performed while Roberto Tempo was visiting Coordinated Science Laboratory at the University of Illinois under CSL Visitor Research Program; the support is gratefully acknowledged.